# Localization in Wireless Sensor Networks Using Quadratic Optimization

Pouya Mollaebrahim Ghari[a], Reza Shahbazian[a], Seyed Ali Ghorashi[a,b]

[a] *Cognitive Telecommunication Research Group, Department of Electrical Engineering, Shahid Beheshti University, Tehran, Iran*
[b] *Cyber Research Centre, Shahid Beheshti University, Tehran, Iran*

p.molaebrahim@mail.sbu.ac.ir, r_shahbazian@sbu.ac.ir and a_ghorashi@sbu.ac.ir

**Abstract** The localization problem in a wireless sensor network is to determine the coordination of sensor nodes using the known positions of some nodes (called anchors) and corresponding noisy distance measurements. There is a variety of different approaches to solve this problem such as semi-definite programming (SDP) based, sum of squares and second order cone programming, and between them, SDP-based approaches have shown good performance. In recent years, the primary SDP approach has been investigated and a variety of approaches are proposed in order to enhance its performance. In SDP approaches, errors in approximating the given distances are minimized as an objective function. It is desirable that the distribution of error in these problems would be a delta distribution, which is practically impossible. Therefore, we may approximate delta distribution by Gaussian distribution with very small variance. In this paper, we define a new objective function which makes the error distribution as similar as possible to a Gaussian distribution with a very small variance. Simulation results show that our proposed method has higher accuracy compared to the traditional SDP approach and other prevalent objective functions which are used such as least squares. Our method is also faster than other popular approaches which try to improve the accuracy of the primary SDP approach.

## 1. Introduction

A sensor network may contain a large set of inexpensive and autonomous sensor nodes that can communicate with their neighbors within a limited radio range. Today, sensor networks have a wide variety of applications in diverse areas such as wildlife monitoring, industrial control and monitoring [1], robotics [2] and medical applications [3]. Determination of the exact location of sensors is a necessity in many of these applications.

One approach to solve this problem is using GPS (Global Positioning System). However, this could be expensive and it may not work very well in some areas [4, 7]. Alternatively, the known locations of some nodes (called anchors) and the distance measurements among neighbor sensors [4, 8, 10] can be used to localize the sensors. In this method, distance measurements are usually corrupted by noise, which may cause errors in target location estimation. In order to minimize this error, sensor network localization can be viewed as an optimization problem, and although it is not convex, some relaxation methods are proposed in the literature to solve this problem.

One of the first proposed convex optimization models for wireless sensor network (WSN) localization is introduced in [11] in which the problem is relaxed to a second order cone programming (SOCP) problem. However, this method needs a large number of anchors on the area boundary to show effective performance as reported in [1]. In [5] the authors introduce a semi-definite programming (SDP) relaxation method for WSN localization, and they report that their proposed method outperforms the method in [11] in terms of accuracy. However, some localization errors happen with rather high variances in this method. In [6], two formulations are developed to deal with noisy corrupted data, and SDP relaxation is used to transform the problem into a convex optimization problem; the first approach is based on maximum likelihood (ML) estimation and the second one is based on solving an SDP feasibility problem with upper and lower bound distance measures [6]. However, ML based methods are very time consuming in comparison with other proposed methods.

The authors in [12] propose a method called smaller SDP (SSDP) to further relax the SDP relaxation. In this method, they relax a single semi-definite matrix cone into a set of small-size semi-definite matrix cones. The two relaxations are proposed in [12] include a node-based relaxation (called NSDP) and an edge-based relaxation named ESDP. Although these proposed relaxations are weaker than the original SDP one, the SSDP is faster compared to others. To find a low rank solution for this problem is investigated in [13] in which a non-convex objective function is proposed. However, the use of this objective function failed to provide more accuracy compared to traditional SDP models. Authors in [2] propose a convex relaxation to solve the localization problem based on maximum likelihood (ML) formulation. The convex relaxation is similar to the SSDP. In [14] the sum of squares (SOS) convex relaxation method is proposed which is very accurate, but it is at the cost of high computation complexity. In [15], authors attempt to address the non-convexity of the problem by introducing a convex objective function but this method is only effective for noise-free measurement cases.

Minimizing the sum of absolute errors as an objective function in the optimization problem may cause large errors [5], and this affects the accuracy of the localization process (this will be shown in section 4). In this paper, we propose a multi-criterion optimization method in WSN localization to reach greater accuracy with not much computational complexity compared to those methods mentioned earlier. Our proposed objective function makes the distribution of error close to Gaussian distribution with a very small variance. This is an approximation of desirable delta distribution for estimation error and therefore the proportion of large errors is effectively reduced.

The remainder of this paper is organized as follows. Section 2 explains the problem formulation. In section 3 the proposed convex relaxations are introduced. The numerical results of the proposed technique and its performance comparisons are presented in section 4, and finally section 5 concludes the paper.

## 2. Problem Formulation

In this section we present the Biswas-Ye model for WSN localization. In this method, the problem is relaxed to a SDP model. We use $I$, $e_i$ and $\mathbf{0}$ to denote the identity matrix, zero column vector except for 1 in position $i$ and the vector of all zeros, respectively. The 2-norm of a vector $x$ denoted as $\|x\|$. Also a positive semi-definite matrix $X$ is represented by $X \succeq 0$. For the sake of simplicity, let the sensor points be placed on a plane. The generalization of this problem for higher dimensions is straightforward. Assume that there are $m$ points with known locations (anchors) $a_k \in \mathcal{R}^2, k=1,\ldots,m$ and $n$ unknown points (sensors) $x_i \in \mathcal{R}^2, i=1,\ldots,n$. $\hat{d}_{ij}$ is the distance measurement between two sensors $x_i$ and $x_j$, and $\hat{d}_{kj}$ is the distance measurement between some sensor $x_j$ and anchor $a_k$.

$$\|x_j - x_i\|^2 = \hat{d}_{ij}^2, \forall (j,i) \in N_s$$
$$\|x_j - a_k\|^2 = \hat{d}_{jk}^2, \forall (j,k) \in N_a$$

where $N_s$ and $N_a$ are the set of pairs $(j,i)$ in which the sensors $i$ and $j$ are connected and the set of pairs $(j,k)$ in which the sensor $j$ is connected to anchor $k$, respectively. Nodes $(i,j) \in N_s$ (or $(j,k) \in N_a$) are called *connected* if the associated measurements $\hat{d}_{ij}$ (or $\hat{d}_{kj}$) do exist for sensor/anchor pairs. Then, the localization problem is to find the location of $n$ sensors by means of given $m$ anchor locations and some distance measurements. We consider a two dimensional localization and then extension to three dimensional is straight forward. Hence, the localization problem can be written as follows:

$$\text{find} \quad X \in \mathcal{R}^{2\times n}, Y \in \mathcal{R}^{n\times n} \qquad (1)$$
$$s.t. (e_i - e_j)^T Y (e_i - e_j) = \hat{d}_{ij}^2, \forall (j,i) \in N_s$$

$$\begin{pmatrix} a_k \\ -e_j \end{pmatrix}^T \begin{pmatrix} I_2 & X \\ X^T & Y \end{pmatrix} \begin{pmatrix} a_k \\ -e_j \end{pmatrix} = \hat{d}_{jk}^2, \forall (j,k) \in N_a$$

where $I_2$ is a 2 by 2 identity matrix, matrix $X = [x_1, \dots, x_n]$, $Y = X^T X$. The problem (1) is non-convex and may relaxed by substituting the constraint $Y = X^T X$ into $Y \succeq X^T X$. This matrix inequality is equivalent to [19]:

$$Z = \begin{pmatrix} I_2 & X \\ X^T & Y \end{pmatrix} \succeq 0 \qquad (2)$$

We define:

$$A_I = \begin{pmatrix} \mathbf{0} & \mathbf{0} & \mathbf{0} \\ 1 & 0 & 1 \\ 0 & 0 & 1 \end{pmatrix}, b_I = \begin{pmatrix} 1 \\ 1 \\ 2 \end{pmatrix} \qquad (3)$$

where **0** in $A_I$ is a zero column vector of dimension $n$ and by introducing slack variables $\alpha^+$ and $\alpha^-$ [5], the problem in Eq. (1) can be written as an standard SDP problem:

$$\min \sum_{(j,i)\in N_s} (\alpha_{ij}^+ + \alpha_{ij}^-) + \sum_{(j,k)\in N_a} (\alpha_{jk}^+ + \alpha_{jk}^-) \qquad (4)$$

$$s.t.\ diag\left(A_I^T Z A_I\right) = b_I$$

$$\begin{pmatrix} e_i - e_j \\ \mathbf{0} \end{pmatrix}^T Z \begin{pmatrix} e_i - e_j \\ \mathbf{0} \end{pmatrix} - \alpha_{ij}^+ + \alpha_{ij}^- = \hat{d}_{ij}^2, \forall (j,i) \in N_s$$

$$\begin{pmatrix} e_i \\ -a_k \end{pmatrix}^T Z \begin{pmatrix} e_i \\ -a_k \end{pmatrix} - \alpha_{jk}^+ + \alpha_{jk}^- = \hat{d}_{jk}^2, \forall (j,k) \in N_a$$

$$\begin{pmatrix} e_i - e_j \\ \mathbf{0} \end{pmatrix}^T Z \begin{pmatrix} e_i - e_j \\ \mathbf{0} \end{pmatrix} \geq r_{ij}^2, \forall (j,i) \notin N_s$$

$$\begin{pmatrix} e_i \\ -a_k \end{pmatrix}^T Z \begin{pmatrix} e_i \\ -a_k \end{pmatrix} \geq r_{jk}^2, \forall (j,k)) \notin N_a$$

$$Z \succeq 0$$

$$\alpha_{ij}^+, \alpha_{ij}^-, \alpha_{jk}^+, \alpha_{jk}^- \geq 0$$

$$for\ i, j = 1, \dots, n, k = 1, \dots, m$$

where $r_{ij}$ or $r_{jk}$ is the radio range in which the associated nodes can communicate. Biswas and Ye [5] neglect the inequalities involving $r_{ij}$ and $r_{jk}$ in order to reduce the problem size; it is shown in [1] that the resulting problem is likely to satisfy most of neglected constraints.

## 3. Proposed quadratic programming approach

In [5] the objective function which is used in the localization problem is the $l_1$-norm of the squared distance error. The authors of [9] reported that if the least-squares penalty function is picked as an objective function, more accuracy than using $l_1$-norm can be achieved. On one hand, it is known that the least-squares objective function put small weights on small squared-distance errors [16]. Therefore, the localization model is comparatively sensitive to large squared-distance errors and reducing the number of large squared-distance errors can improve the performance of the model in average. On the other hand, it is ideal that the squared-distance error has a delta shaped PDF (probability density function) which is impractical. In this paper, we aim to introduce a new objective function which tries to estimate the squared-distance error's PDF as a Gaussian distribution with a small variance, in order to avoid large squared-distance errors and simultaneously maintain proportionally small amounts of squared-distance errors. To start, we consider the following theorem [17]:

**Theorem 1.**For any random variable $X$ and its estimate $\hat{X}$,

$$E(X - \hat{X})^2 \geq \frac{1}{2\pi e} e^{2h(X)} \qquad (5)$$

where $h(X)$ is the differential entropy of random variable $X$ and $\hat{X}$ is the mean value of $X$. Equality happens in Eq. (5) if and only if $X$ is Gaussian.

Theorem 1 states that the variance of desirable estimation is zero and it is a delta function. It is impossible to have no estimation error in practical scenarios. Therefore, we reduce the variance of the estimation error which is more likely to reduce the relative entropy between the PDF of squared-distance errors and Gaussian distribution. Since the Gaussian distribution is characterized by its first and second statistics, we aim to approximate these statistics in order to make error distribution close to Gaussian distribution with as small a standard deviation as possible. Thus, involving the variance of squared-distance error in the optimization problem can reduce the difference between the distribution of squared-distance error and Gaussian distribution with a small standard deviation.

As mentioned earlier in section 1, we consider two criteria for optimization in order to obtain higher accuracy. The first criterion is the mean of the squared-distance errors i.e. $\alpha_{ij}^{+}-\alpha_{ij}^{-}$ and $\alpha_{jk}^{+}-\alpha_{jk}^{-}$. Obviously, we aim to minimize this criterion in order to make the mean of squared-distance errors as close as possible to the mean of associated Gaussian distribution which is zero. The second criterion is the variance of the squared-distance errors, in order to make the PDF of the squared distance errors as close as possible to delta distribution. The delta distribution is approximated by Gaussian distribution and we consider the variance of squared-distance errors as a second criterion for the optimization.

By interpreting the localization problem as a bi-criterion optimization, we have to minimize a vector $f_0$, which consists of two scalar functions $F_1$ and $F_2$ where $F_1$ is the mean and $F_2$ is the variance of the squared distance error as follows:

$$F_1=\frac{1}{v}\left(\sum_{(j,i)\in N_s}\left|\alpha_{ij}^{+}-\alpha_{ij}^{-}\right|+\sum_{(j,k)\in N_a}\left|\alpha_{jk}^{+}-\alpha_{jk}^{-}\right|\right) \quad (6)$$

$$F_2=\frac{1}{v}\left(\sum_{(j,i)\in N_s}(\left|\alpha_{ij}^{+}-\alpha_{ij}^{-}\right|-F_1)^2+\sum_{(j,k)\in N_a}(\left|\alpha_{jk}^{+}-\alpha_{jk}^{-}\right|-F_1)^2\right) \quad (7)$$

where $v$ is the number of distances between connected nodes. The objective function $f_0$ must be transformed into a scalar one by vector $\lambda$ ( $\lambda\succ 0$ ) in a way that the optimization problem can be solved. Therefore, we may write the objective function as:

$$F_1+\gamma F_2 \quad (8)$$

This objective function is non-convex. Therefore, we relax it to a convex form as follows:

$$\gamma\sum_{(j,i)\in N_s\bigcup N_a}(\left|\alpha_{ij}^{+}-\alpha_{ij}^{-}\right|+\frac{1}{2\gamma})^2-2F_1\gamma\sum_{(j,i)\in N_s\bigcup N_a}\left|\alpha_{ij}^{+}-\alpha_{ij}^{-}\right|+\gamma vF_1^2-\frac{v}{4\gamma} \quad (9)$$

Assume that $\alpha\in\mathcal{R}^v$ is a vector with components of $\alpha_{ij}^{+}-\alpha_{ij}^{-}$ where $(i,j)\in N_s\bigcup N_a$. Then, we approximate $F_1$ by constant $m'$ that satisfies the following constraints. Let *l* be a constant value which satisfies the following:

$$m'\le\frac{1}{v}\sum_{(j,i)\in N_s\bigcup N_a}\alpha_{ij}\le\frac{1}{v}\sum_{(j,i)\in N_s\bigcup N_a}\left|\alpha_{ij}\right| \quad (10)$$

$$m'\le\sqrt{1-\frac{1}{v}}m \quad (11)$$

$$m=\frac{1}{v}\sum_{(j,i)\in N_s\bigcup N_a}\alpha_{ij} \quad (12)$$

$$F_1^2-m^2\le l$$

We relax the objective function by using the Jensen's inequality and write the relaxed objective function as follows:

$$\frac{\gamma}{v}\left\|\alpha+\frac{1}{2\gamma}\right\|^2+\frac{\gamma}{v^2}\left\|\alpha-\sum_{(j,i)\in N_s\bigcup N_a}\alpha_{ij}\right\|^2-\frac{1}{4\gamma}-\gamma m'^2v+v\gamma l \quad (13)$$

In Eq. (13) we have substituted $\frac{1}{v}\sum_{(j,i)\in N_s\bigcup N_a}\alpha_{ij}$ with $F_1$. Since the objective function is positive, this tightening preserves convexity. Now we can conclude that:

$$m'\le\sqrt{1-\frac{1}{v}}m\rightarrow m^2\left(v-1\right)-vm'^2\ge 0\rightarrow\gamma v(m^2v-m'^2v-m^2)\ge 0 \quad (14)$$

If $\gamma$ satisfies the constraint $\gamma \geq 1/2m$ along with other constraints in the problem, we prove that the relaxed objective function is an upper bound of the non-convex objective function. Knowing the convex upper bound, we may solve the optimization problem. Therefore we have:

$$(2m\gamma-1)v \sum_{(j,i)\in N_s \bigcup N_a} (|\alpha_{ij}| - \alpha_{ij}) + \gamma \sum_{(j,i)\in N_s \bigcup N_a} \alpha_{ij}^2 + \gamma m^2 v^2 - \gamma m'^2 v^2 - \gamma m^2 v \geq 0 \qquad (15)$$

We may write the following statement:

$$\frac{\gamma}{v}\left\|\alpha + \frac{1}{2\gamma}\right\|^2 + \frac{\gamma}{v^2}\left\|\alpha - \sum_{(j,i)\in N_s \bigcup N_a} \alpha_{ij}\right\|^2 - \frac{1}{4\gamma} - \gamma m'^2 v + v\gamma l \geq F_1 + \gamma F_2 \qquad (16)$$

After eliminating constant terms and relaxing the objective function, we can write the objective function as follows:

$$v\left\|\alpha + \frac{1}{2\gamma}\right\|^2 + \left\|\alpha - \sum_{(j,i)\in N_s \bigcup N_a} \alpha_{ij}\right\|^2 \qquad (17)$$

We take $\gamma$ as $\frac{2v}{\sum_{(j,i)\in N_s \bigcup N_a} \alpha_{ij}}$ in Eq. (17) to satisfy all the constraints, therefore:

$$v\left\|\alpha + \frac{1}{v}\sum_{(j,i)\in N_s \bigcup N_a} \alpha_{ij}\right\|^2 + \left\|\alpha - \sum_{(j,i)\in N_s \bigcup N_a} \alpha_{ij}\right\|^2 \qquad (18)$$

Since the regularization coefficient $\gamma$ must be positive and in order to guarantee that the relaxed objective function is the upper bound of the original objective function, we have to add the following constraint in the SDP model:

$$\sum_{(j,i)\in N_s \bigcup N_a} \alpha_{ij} \geq 0 \qquad (19)$$

On one hand, adding Eq. (19) as a constraint in the optimization model increases the problem size. On the other hand, the resulting model with neglecting this constraint is likely to satisfy both positive $\gamma$ and the objective function upper bound. Therefore, we may neglect this constraint.

It is straight forward to prove that the function in Eq. (18) is convex. Using the relaxation Biswas and Ye [5] applied for sensor network localization, we may write the optimization problem as Eq. (20) where $w = mean(\alpha)$ in the following model:

$$\begin{aligned}
&\min vt_1 + t_2 \qquad (20)\\
&s.t.\ diag\left(A_I^T Z A_I\right) = b_I\\
&\begin{pmatrix} e_i - e_j \\ \mathbf{0} \end{pmatrix}^T Z \begin{pmatrix} e_i - e_j \\ \mathbf{0} \end{pmatrix} - \alpha_{ij}^+ + \alpha_{ij}^- = \hat{d}_{ij}^2, \forall (j,i) \in N_s\\
&\begin{pmatrix} e_i \\ -a_k \end{pmatrix}^T Z \begin{pmatrix} e_i \\ -a_k \end{pmatrix} - \alpha_{jk}^+ + \alpha_{jk}^- = \hat{d}_{jk}^2, \forall (j,k) \in N_a\\
&\|\alpha + w\| \leq s_1\\
&\|\alpha - vw\| \leq s_2\\
&\begin{pmatrix} 1 & s_1 \\ s_1 & t_1 \end{pmatrix} \succeq 0 \qquad (20.a)\\
&\begin{pmatrix} 1 & s_2 \\ s_2 & t_2 \end{pmatrix} \succeq 0 \qquad (20.b)\\
&Z \succeq 0\\
&\alpha_{ij}^+, \alpha_{ij}^-, \alpha_{jk}^+, \alpha_{jk}^- \geq 0\\
&for\ i, j = 1, \ldots, n, k = 1, \ldots, m
\end{aligned}$$

Eq. (20.a) and Eq. (20.b) are the SOC constraints that are equivalent to quadratic constraints.

## 4. Simulation Results

In this section, we display the simulation results and compare our proposed quadratic programming (QP) with other models. We use MATLAB R2010b and also SDPT3 solver in CVX software [18] in order to solve the problems. The CPU time which is taken for localizing is extracted from

cvx_cputime [18]. We consider 2-dimentional problems and use the benchmark test10-500 which is available online at http://www.stanford.edu/~yyye. The position error for sensor $i$ is defined as $\|x_i - \hat{x}_i\|$ where $x_i$ is the true position for sensor $i$ and $\hat{x}_i$ is the estimated location for this sensor. *PE* shows the position error for the evaluated network and is defined as:

$$PE = \frac{1}{n}\sum_{i=1}^{n}\|x_i - \hat{x}_i\|$$

At first we seek to find appropriate regularization coefficient in terms of position error and survey the effect of our proposed QP model on relative entropy.

### 4.1. Regularization Coefficient

The section mentions the evaluation results of the regularization coefficient effect on the position error. We consider networks which include 80 sensors and 5 anchors. The distance measurements are corrupted by an additive Gaussian noise with standard deviation of 0.05 and radio range is set to 0.25. The experiments are performed over 50 networks. As depicted in Fig. 1, when the regularization coefficient becomes larger, the average position error converges to a limit. As can be seen in Fig. 1, if the regularization coefficient is large enough, the variation of it doesn't cause a rapid change in average position error. Therefore, we can solve the problem regardless of regularization coefficient as long as the regularization coefficient is large enough. We try to approximate this specific amount.

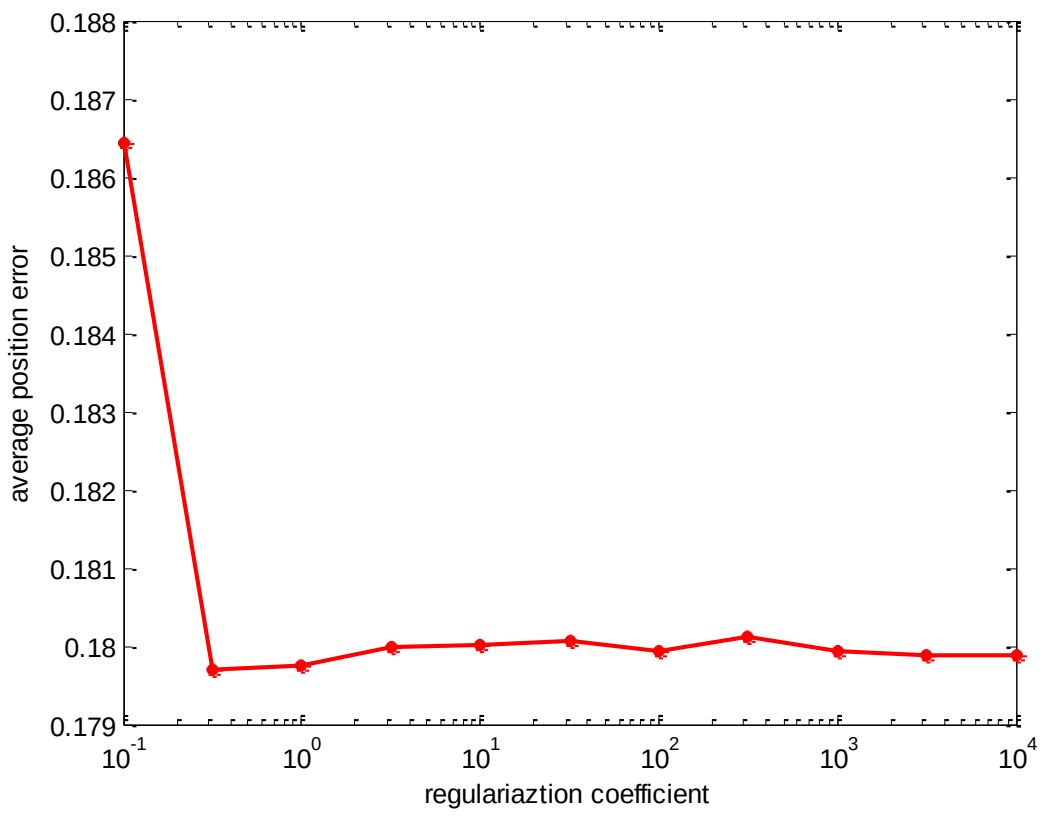

**Fig. 1.** Average of position error as a function of regularization coefficient variations.

### 4.2. Relative Entropy

Relative entropy $D$ is the criterion that measures the difference between two distributions $p$, $q$ and is defined as follows:

$$D(p\,||\,q) = \sum_{x\in\chi} p(x)\log\frac{p(x)}{q(x)}$$

In this section, we want to prove that our model reduces the relative entropy between distribution of squared-distance error and Gaussian distribution. We quantize the probability distribution of squared-distance error and divide real numbers into intervals $I_0$, $I_1$, ... with a specific length and then calculate the probability associated to interval $I_i, i \in \{0,1,\ldots\}$ as follow:

$$p_{I_i} = v_{I_i}/v$$

where $v_{I_i}$ is the number of squared-distance errors lie in $I_i$ and $v$ is the total number of squared-distance errors.

Now, we study the effect of our proposed QP approach on the amount of mentioned relative entropy. We set the length of intervals to 0.0049 and choose 0.008 as the standard deviation of Gaussian distribution. Networks consist of 80 sensors and 5 anchors. The distance measurements between connected nodes are corrupted by additive Gaussian noise with standard deviation 0.05 and similar to prior experiment, the radio range is set to 0.25. We use 50 networks for operating the experiment.

The relative entropy is measured for the proposed QP, Biswas-Ye and LS. The average of results are depicted in Table 1 and show that our proposed QP approach reduces the relative entropy dramatically in comparison with the approach in [5] and it has also a lower relative entropy in comparison with LS.

**Table 1 .** Relative entropy of squared-distance errors in bits for Biswas-Ye [5], LS [9] and our proposed QP.

| **Approach** | **Relative entropy (bits)** |
|---|---|
| Biswas-Ye [5] | 1.1232 |
| LS [9] | 0.2229 |
| Proposed QP | 0.2089 |

Figs. 2,3,4 depict the distribution of squared-distance errors for Biswas-Ye, LS and our proposed QP approaches, respectively for one of the 50 simulated networks and illustrate that the distribution of squared-distance errors for proposed QP is more similar to Gaussian distribution. The proportion of squared-distance errors in the interval [0.00245, 0.0073) for the proposed QP is 0.3598 and for LS is 0.3482. The proportion of squared-distance errors in the interval [-0.00245, 0.00245) is 0.2244 for both methods in this case. These results along with comparisons of Fig. 3 and Fig. 4 show that the density of squared distance error in the vicinity of 0 in proposed QP is more than what in LS. As can be seen in Figs. 2, 3 and 4 the density of large squared-distance errors in Biswas-Ye is more than two other approaches. For more illustration consider one of the simulated networks. The proportion of squared-distance errors which are larger than 0.022 for Biswas-Ye approach is 0.0658, for LS is 0.0271 and in our proposed QP is 0.0348. In addition the maximum amount of squared distance error in Biswas-Ye approach is more than LS and our QP model. Large squared-distance errors can degrade the performance of localization effectively.

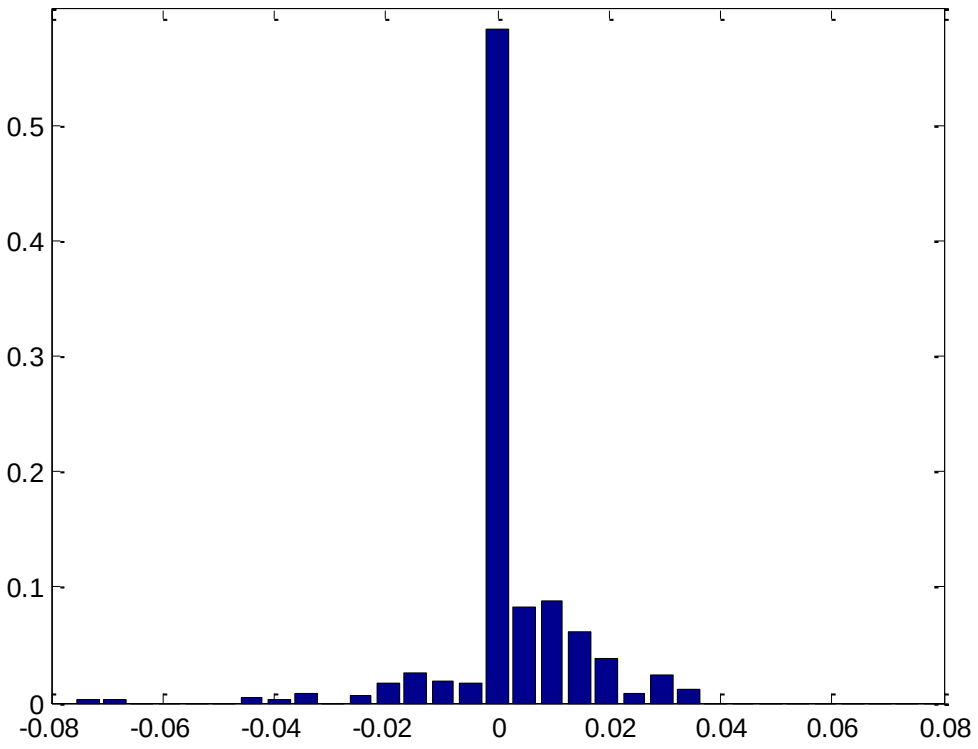

**Fig.2.** Distribution of squared-distance errors in Biswas-Ye approach [5].

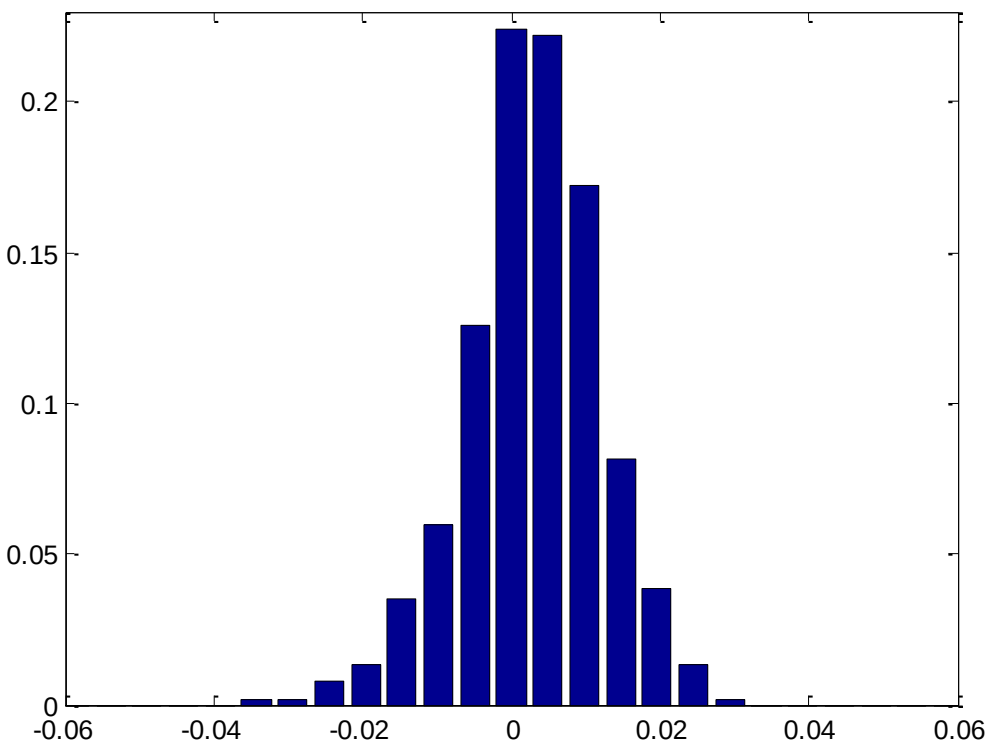

**Fig.3.** Distribution of squared-distance errors in LS approach [9].

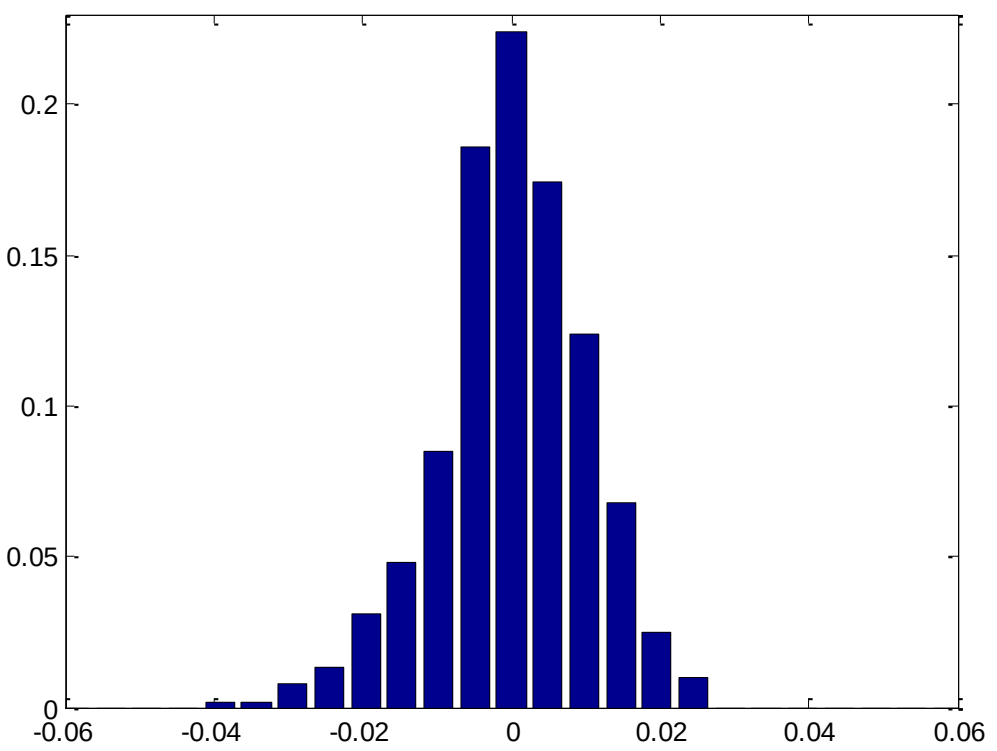

**Fig.4.** Distribution of squared-distance errors in proposed QP approach.

### 4.3. Noise Standard Deviation

In this section, we consider the performance of our proposed QP approach in comparison with Biswas-Ye, EML, ML and LS by increasing the noise value. We evaluate our proposed method and other introduced methods over 50 networks. As depicted in Fig. 5, the performance of our proposed QP shows higher accuracy in comparison with Biswas-Ye, LS and EML in all noise standard deviation intervals. As can be seen in Fig. 5, when the standard deviation of noise becomes larger, the difference between the accuracy of proposed QP and ML reduces. However, as will be shown in 4.5 our QP method has lower complexity in comparison with ML.

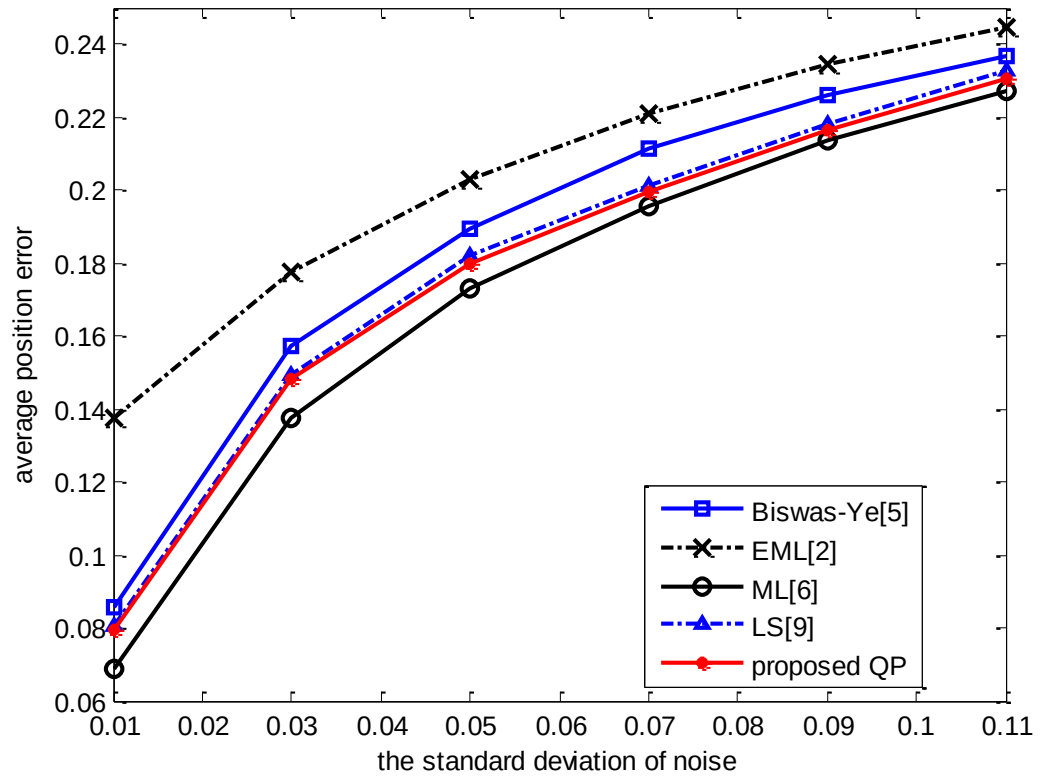

**Fig.5.** Performance comparison between the approaches for various standard deviations of noise.

### 4.4. Radio Range

In this section, we report the performance evaluation results of Biswas-Ye, EML, ML LS and proposed QP by changing radio range. We set the number of sensors to 80 and anchors to 5, the standard deviation of noise is assumed to be 0.05 and we perform the experiment over 50 networks.

In Fig. 6 we compare the average error of our proposed QP approach with Biswas-Ye [5], EML, ML and LS. As depicted in Fig. 6, increasing the radio range does not affect the EML properly. It can be seen that our proposed QP approach have lower position error compared with those introduced in [2, 5]. As we can see in Fig. 6, the performance of our proposed QP is better than the one of LS [9]. However, we will discuss in section 4.5 that the proposed QP has less computational complexity.

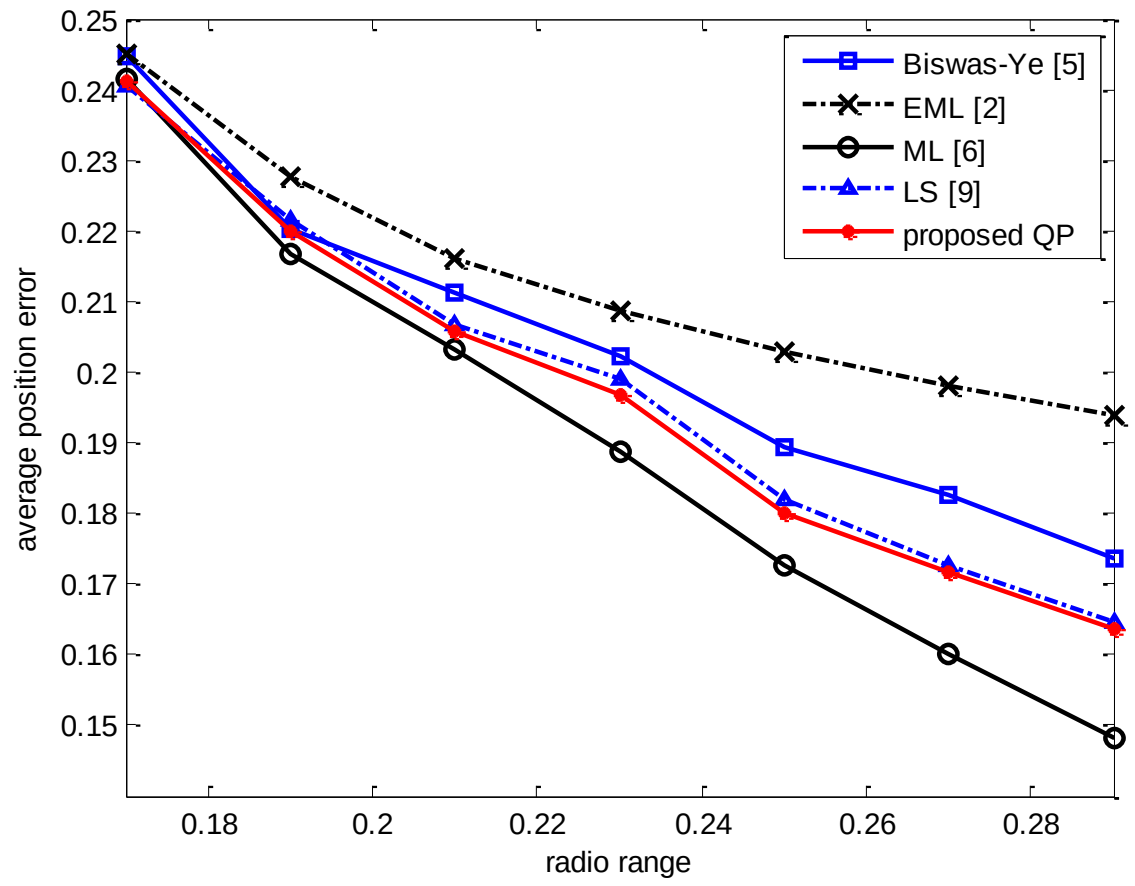

**Fig.6.** Average position error of the approaches in various radio ranges.

### 4.5. Number of Sensors and Solution Time

In this section we would evaluate the effect of number of sensors on problem solution time. In this simulation, the radio range is set to 0.3 and distance measurements are corrupted by additive Gaussian noise with standard deviation of 0.05 and experiments are performed over 50 networks. By increasing the number of sensors, the number of constraints and consequently solution time increases. In practical scenarios as the number of sensors increases, the required radio range and the number of constraints scales with $O(n)$. Therefore, the number of operations is typically $O(n^3)$ [6]. Fig. 7 depicts the effect of increasing the number of sensors on solution time. Fig. 7 illustrates that the solution time of our proposed QP method is lower than the solution time of LS, ML and EML. Although its solution time is more than the solution time of Biswas-Ye which is acceptable by showing higher accuracy as depicted in Figs. 5 and 6. As it can be observed in Fig. 7, the solution time of Biswas-Ye changes approximately linearly by increasing the number of sensors while solution time of ML and EML methods change sharper compared to LS and our proposed QP and we conclude that they have a higher level of complexity.

However, the solution time can be decreased by limiting the number of connected sensors or equivalently limiting the connectivity. In Fig. 8 we study the effect of limiting the connectivity on the solution time. We set the connectivity to 7 i.e. the maximum number of connected sensors for each sensor is 7. As can be seen, EML shows better performance in comparison with its performance in Fig.7. Our proposed QP has lower solution time compared to ML, LS and EML. In addition, Fig.8 illustrates that by applying the limitation on connectivity, the solution time decreases considerably. However, limiting the connectivity degrades the performance of localization in terms of accuracy. For example for 50 networks consist of 100 sensors, the average position error of ML with considering limitation in connectivity is 0.1727 and in the case without limiting connectivity, it is 0.1332.

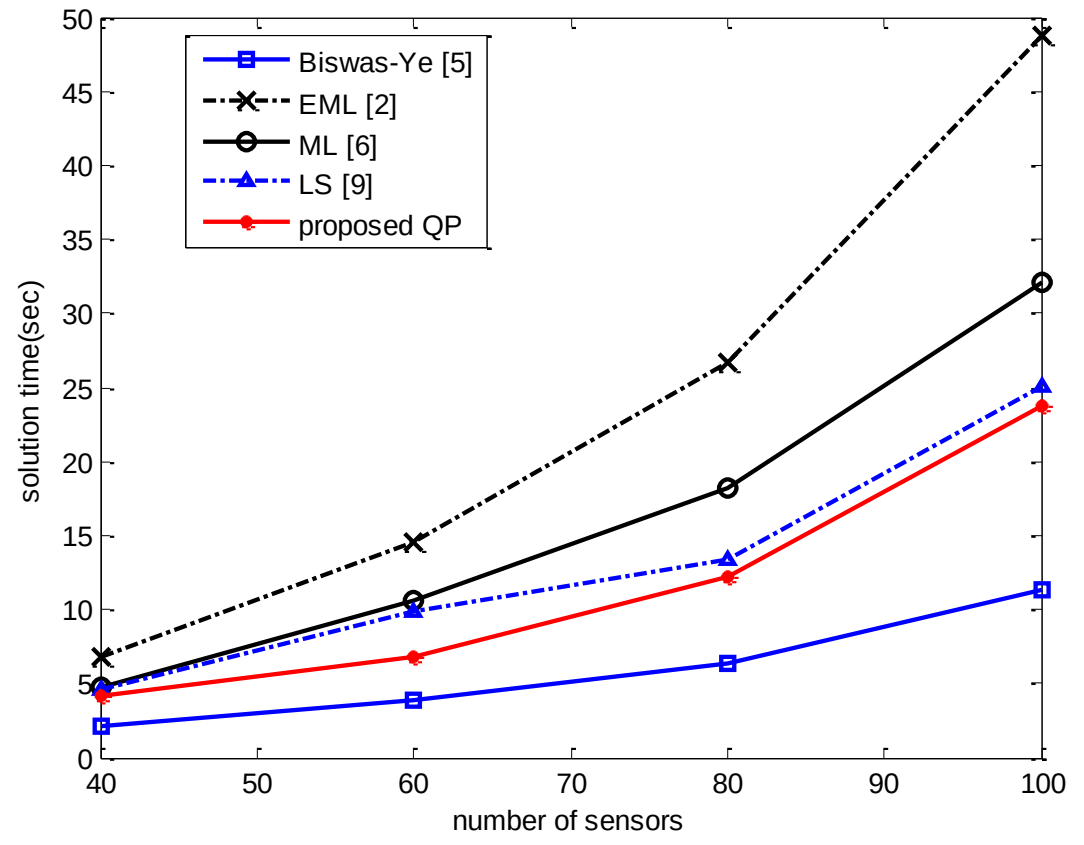

**Fig.7.** The effect of number of sensors on solution time (sec).

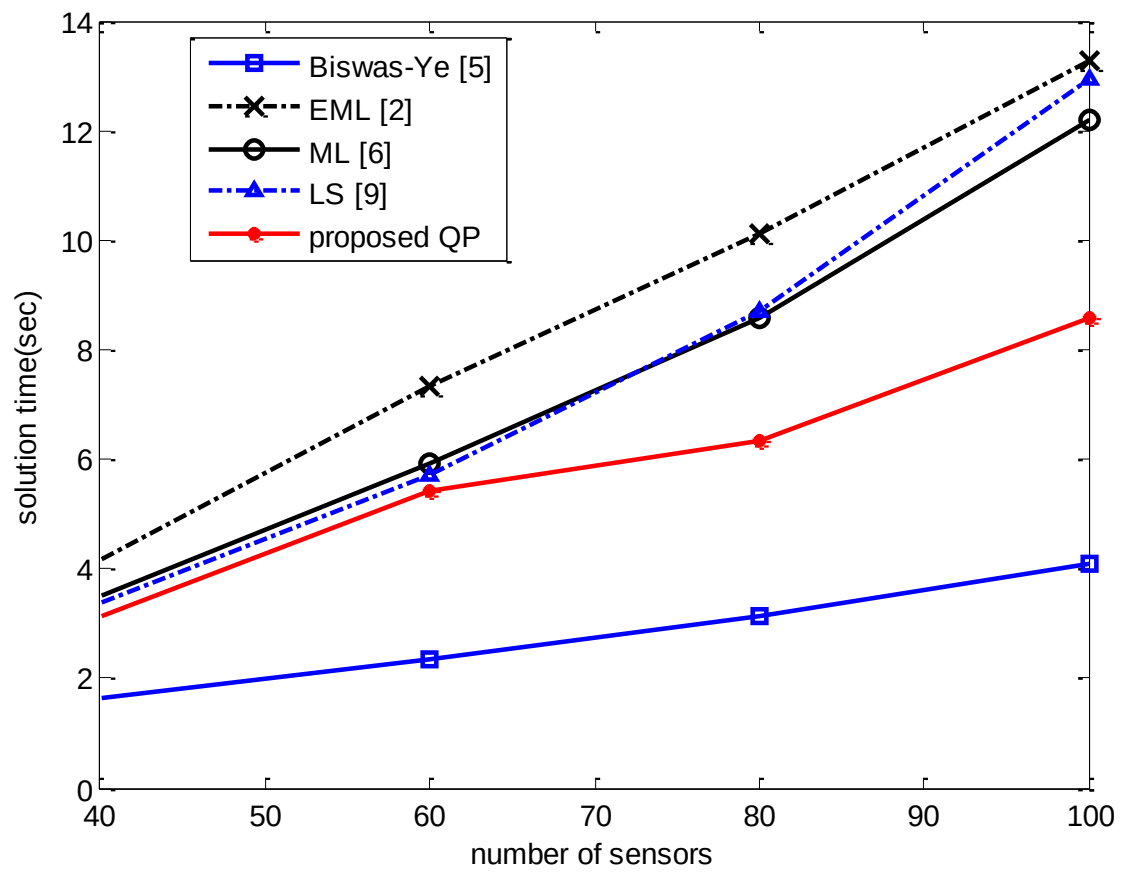

**Fig.8.** The effect of number of sensors on solution time (sec) with limiting the connectivity to 7.

Fig. 9 depicts the coordinates of sensor's location estimation for Biswas-Ye, EML, ML, LS and proposed QP approach along with exact locations. The network consists of 20 sensors and 5 anchors, radio range is fixed to 0.4 and the standard deviation of additive Gaussian noise is 0.01. As can be seen in Fig. 9, EML has worst performance among all other approaches, LS and proposed QP has the same performance in most cases.

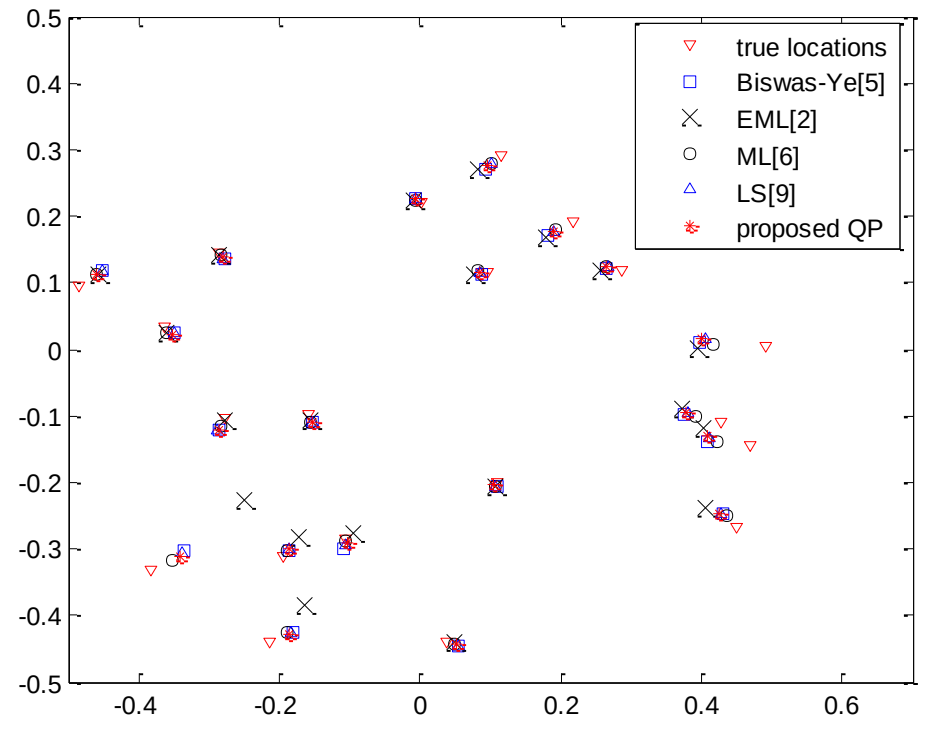

**Fig.9.** The coordinates which localize with Biswas-Ye [5], EML[2], ML[6], LS[9] and proposed QP along with their exact locations.

## 5. Conclusion

In this paper we consider the variance of error in order that we can define a new objective function. Applying this new objective function makes the error distribution as similar to a Gaussian distribution as possible with a very small variance. To do this we introduce a bi-criterion objective function which is non-convex but we relax it to a convex form in order to solve the problem. Our proposed objective function reduces proportion of large errors effectively and also increases the proportion of errors in vicinity of zero compared to least squares function. Simulation results show that our proposed QP method improves the accuracy in position estimation compared to other prevalent objective functions such as least squares one. Computational complexity is also less than ML-based methods like EML, ML and a bit more compared to Biswas-Ye, which is worth it due to its performance and accuracy in the localization problem. Therefore, simulation results show that our proposed QP is faster than those that attempt to enhance the accuracy of the SDP approaches.